\documentclass[10pt]{article}
\usepackage{amsmath}
\usepackage{amsfonts}
\usepackage{amssymb}
\usepackage{comment}
\usepackage{graphicx,color}
\usepackage{authblk}
\long\def\proof#1{\removelastskip\vskip\baselineskip\relax\noindent{\it
Proof\if!#1!\else\ \ignorespaces#1\fi.\ }\ignorespaces}

\newcommand{\Q}{{\mathbb Q}}
\newcommand{\Z}{{\mathbb Z}}

\newcommand{\C}{{\mathbb C}}

\newcommand{\z}{\zeta}

\DeclareMathOperator{\asin}{asin}

\newcommand{\squareforqed}{\hbox{\rlap{$\sqcap$}$\sqcup$}}
\newcommand{\qed}{\ifmmode\squareforqed\else{\unskip\nobreak\hfil
\penalty50\hskip1em\null\nobreak\hfil\squareforqed
\parfillskip=0pt\finalhyphendemerits=0\endgraf}\fi}

\newtheorem{theorem}{Theorem}[section]

\newtheorem{proposition}[theorem]{Proposition}

\newtheorem{definition}[theorem]{Definition}
\newtheorem{conjecture}[theorem]{Conjecture}

\newcommand{\litem}{\par\noindent\dimen0=\parindent%
    \advance\dimen0 by-4pt
               \hangindent=\dimen0\ltextindent}

\newcommand{\ltextindent}[1]{\hbox to \hangindent{#1\hss}\ignorespaces}
\newcommand{\ltextjndent}[1]{\hbox to \hangindent{#1\hss}\ignorespaces\kern-1ex}

\definecolor{lightgrey}{rgb}{0.8, 0.84, 0.8}

\begin{document}
\pagestyle{plain}

\title{Parametric Continued Fractions for $\pi^2$, $\z(3)$, and other Constants}
\author{Henri Cohen}
\affil{Universit\'e de Bordeaux, LFANT,\\
  IMB, U.M.R. 5251 du C.N.R.S., 351 Cours de la Lib\'eration\\
33405 Talence Cedex, FRANCE}

\maketitle

\begin{abstract}
  We give an extensive list of parametrized \emph{families} of polynomial
  continued fractions of \emph{smallest possible degrees} for $\pi^2$ and
  $\z(3)$, and mention similar results for other constants.
\end{abstract}

\smallskip

\section{Introduction}

One can find in the literature a huge number of continued fractions
(abbreviated CFs) for interesting constants and functions. The purpose of
this paper is to give possibly \emph{experimental} results showing that there
often exist \emph{families} of such CFs, usually depending on two, three,
and even four parameters. In many cases these experimental results can be
\emph{proved}, usually because either $p_n/q_n$ or $p_n/q_n-p_{n-1}/q_{n-1}$
have an explicit expression. In other cases, I have not bothered to try
to prove the experimental results, although it may not be difficult.
Thus, when I speak of results, they should be considered experimental, but
with considerable evidence. I would appreciate receiving any additional
data of the same kind.

\smallskip

Notation: a CF will be written as a pair of sequences $(a_n,b_n)_{n\ge0}$
representing the CF $a_0+b_0/(a_1+b_1/(a_2+b_2/(a_3+\cdots)))$ (warning:
note this precise convention, other authors have different notation).
All the CFs that we will encounter will be \emph{convergent}: if as usual
$p_n/q_n=a_0+b_0/(a_1+\cdots+b_{n-1}/a_n)$, this means that
$z=\lim_{n\to\infty}p_n/q_n$ exists, and when we speak of the speed of
convergence, we mean the behavior of $z-p_n/q_n$.

Note that we always assume that $b_n\ne0$ for all $n$, but we allow a
finite number of $a_n$ to vanish.

\begin{definition}
\begin{enumerate}\item
  Let $(a_n)_{n\ge0}$ be a sequence (always assumed to be of rational numbers,
  and usually of integers). We say that it is ultimately polynomial if
  there exists a polynomial denoted $a(x)$ (by abuse of notation) with
  rational coefficients such that $a_n=a(n)$ and $a_n\ne0$ for $n$
  sufficiently large.
\item We say that the CF $(a_n,b_n)_{n\ge0}$ is of \emph{polynomial type}
  if the sequences $a_n$ and $b_n$ are ultimately polynomial.
\end{enumerate}
\end{definition}

As a practical notation, if $a_n=a(n)$ for $n\ge n_0$, we will denote the
sequence $(a_n)_{n\ge0}$ by the $(n_0+1)$-component vector
$A=(a_0,a_1,\dotsc,a_{n_0-1},a(n))$, where the last entry is a polynomial.
Of course this is not unique: we also have for instance
$A=(a_0,a_1,\dotsc,a_{n_0-1},a_{n_0},a(n))$ with $a_{n_0}=a(n_0)$.

Thus a CF will be denoted by a pair of such vectors $A$ and $B$ representing
$(a_n)$ and $(b_n)$, and if the CF converges to some number $z$, we will
write $z=(a_n,b_n)$.

\begin{definition} We recall that two CFs $(a_n,b_n)$ and $(a'_n,b'_n)$ are
  said to be \emph{equivalent} if there exists an integer $k$ such that
  $(a'_n,b'_n)=(a_{n+k},b_{n+k})$ for $n$ sufficiently large, i.e., if
  the \emph{tails} of the CFs are ultimately the same.
\end{definition}

Since we assume everything rational, if we set $z=(a_n,b_n)$ and
$z'=(a'_n,b'_n)$, the CFs for $z$ and $z'$ are equivalent if and only if
there exists a M\"obius transformation with integer coefficients sending $z$
to $z'$, in other words $z'=(Az+B)/(Cz+D)$ for integers $A$, $B$, $C$, $D$
with $AD-BC\ne0$.

We can now give our final definition:

\begin{definition} We say that a pair of polynomials $(A(x),B(x))$ with
  rational coefficients \emph{represents} a constant $z$ if there exists
  a CF $(a_n,b_n)$ converging to $z$ such that $a_n=A(n)$ and $b_n=B(n)$
  for $n$ sufficiently large, and in that case we simply write by abuse
  of notation $z=(A(n),B(n))$.\end{definition}

Because of what we have said above, this is almost, but not quite, equivalent
to the fact that the CF $(A(n),B(n))_{n\ge0}$ converges to a M\"obius
transform of $z$ (always assumed with integer coefficients): indeed, some
values of $A(n)$ and/or $B(n)$ may vanish, which kills this property. But
if we avoid these vanishing values by changing the corresponding coefficients
of the CF to any nonzero value, this is indeed true.

\smallskip

Our goal is, given some interesting constant $z$, to find all (or at least
many) pairs $(A,B)$ representing $z$, subject to suitable restrictions.
In particular, recall that if $(a_n,b_n)_{n\ge0}$ is a CF converging to $z$
and $r_n$ is any sequence with $r_0=1$ and $r_n\ne0$ for all $n$, the CF
$(a_nr_n,b_nr_nr_{n+1})_{n\ge0}$ also converges to $z$. We will thus
restrict to polynomials $(A,B)$ of smallest possible degrees.

\smallskip

To make this more precise, we introduce the following definitions:

\begin{definition} We say that $z\in\C$ is a \emph{rational period of
    degree $k$} if it is equal to the sum of a convergent series of the form
  $\sum_{n\ge1}\chi(n)f(n)$, where $\chi$ is a periodic arithmetic function
  taking \emph{rational} values, and $f\in\Q(x)$ is a rational function with
  rational coefficients whose denominator is of degree $k$.
\end{definition}

Note that this definition is \emph{not} compatible with the definition
of periods as given by Kontsevitch--Zagier \cite{Kon-Zag}. Note also that
a $\Q$-linear combination of rational periods of degree $\le k$ is again
a rational period of degree $\le k$.

\smallskip

Examples of periods of degree 1: $\pi$, $\pi/\sqrt{d}$, $\log(2)$,
$\sqrt{2}\log(1+\sqrt{2})$, $\sqrt{3}\log(2+\sqrt{3})$,
$\sqrt{5}\log((1+\sqrt{5})/2)$.

Examples of periods of degree 2: $\pi^2$, $\pi^2/\sqrt{d}$,
$G=L(\chi_{-4},2)$ (Catalan's constant), $G_3=L(\chi_{-3},2)$.

Examples of periods of degree 3: $\pi^3$, $\pi^3/\sqrt{d}$, $\z(3)$.

\smallskip

By using Euler's transformation of series it is clear that a period has a
CF expansion where $a(n)$ and $b(n)$ are polynomials with
rational coefficients for $n$ large. However this is usually uninteresting,
or at least not any more interesting than the series itself. Thus, really
interesting CFs are those for which the degrees of the polynomials are
reasonably small, specifically $\deg(a(n))\le k$ and $\deg(b(n))\le 2k$ for
$n$ sufficiently large (we will say that it has \emph{bidegree} $\le(k,2k)$).
Trivial examples: $S^+=\sum_{n\ge1}1/P(n)$ and
$S^-=\sum_{n\ge1}(-1)^{n-1}/P(n)$ with $\deg(P)=k\ge2$ have the CFs
$S^+=((0,P(n)+P(n-1)),(1,-P(n)^2))$ and $S^-=((0,P(n)-P(n-1)),(1,P(n)^2))$
by Euler, of the required degrees.

Note that if we apply Euler's transformation to
$S=\sum_{n\ge1}(n+1)/n^3=\z(3)+\z(2)$ which is a period of degree $3$,
we obtain a CF of bidegree $(4,8)$. However, there does exist a CF of
bidegree $(3,6)$, for instance we have the (easily proved) CF:
$$\z(3)+\z(2)=((3,2n^3+2n^2+3n+1),(-1,-n(n+1)^5))\;.$$

\smallskip

For most of this paper, we will focus on the cases $z=\pi^2$ (or
equivalently $z=\z(2)=\pi^2/6$) and $z=\z(3)$,
and attempt to classify \emph{all} representations of $\pi^2$ with bidegree
$\le(2,4)$, and of $\z(3)$ with bidegree $\le(3,6)$.
A completely analogous study can be made for other periods such
as those mentioned above, we will briefly mention this in the last section.

Note that, contrary to \cite{Coh}, we are interested in CFs for ``pure''
periods such as those mentioned above, and not in linear combinations,
which is considerably more restrictive since one needs to get rid of
``parasitic'' terms (for instance $\log(2)$ for $\pi^2$, or $\log(2)$
and $\z(2)$ for $\z(3)$, see below).

\medskip

\section{CFs for $\pi^2$: Convergence in $1/n^P$, $P\in\Z_{\ge1}$}

As far as the author is aware, up to equivalence there exist exactly nine
types of convergence of polynomial continued fractions for $\pi^2$ with
bidegree $\le(2,4)$: two families with convergence in $1/n^P$ with $P\in\Z$,
one with convergence in $1/n^P$ with $P\in1/2+\Z$, one with convergence
in $(-1)^n/n^P$, one with convergence $1/(E^nn^P)$ for each of $E=4/3$, $2$,
$4$, and $-8$, and finally the unique CF of Ap\'ery with
convergence in $1/E^n$ for $E=-((1+\sqrt{5})/2)^{10}$. We consider
each of these families in turn.

\medskip

\subsection{First Family}\label{sec:first}

\medskip

Using Euler's transformation of series into CFs applied to
$\pi^2/6=\sum_{n\ge1}1/n^2$, we see that $\pi^2/6$ is given by the CF
$((0,2n^2-2n+1),(1,-n^4))$, with speed of convergence identical to that of
the series, in $1/n$, which we write in abbreviated form as explained above
where we do not give the initial values of $a_n$ and $b_n$
$$\pi^2=(2n^2-2n+1,-n^4)\;.$$
Remark: even though this notation means that there can be a number of
initial values of $a_n$ and $b_n$ which are not given by the general
polynomial formulas, in most of the examples that we will see in this
paper, only $a_0$ and $b_0$ would need to be specified (for instance
here $a_0=0$ and $b_0=1$), the other $a_n$ and $b_n$ for $n\ge1$ are given
by the polynomials.

\smallskip

If we apply the Bauer--Muir transformation, we check that for any nonnegative
integer $k\ge0$ we also have $\pi^2=(2n^2-2n+k^2+k+1,-n^4)$.

This naturally leads us to \emph{search numerically} for other representations
of the form $\pi^2=(A(n),B(n))$, with $A(n)=2n^2+a_1n+a_0$ and $B(n)$ of
degree $4$
with leading coefficient $-1$. We speed up the search by assuming that $B$
is a product of linear factors with integral coefficients (you are welcome to
try other types of polynomials, but I am pretty sure that you will not find
others), and since $(A(n+j),B(n+j))$ is equivalent to $(A(n),B(n))$, we can
assume that $B(n)=-n(n+u_1)(n+u_2)(n+u_3)$ with $0\le u_1\le u_2\le u_3$.

We find a \emph{very large} number of such representations, from which one
can easily extract one and two-parameter families, less easily three-parameter
families, and finally with considerable difficulty, a single four-parameter
family to rule them all:
$\pi^2=(A(n),B(n))$ with
\begin{align*}
  A(n)&=2n^2+(3u+v+w-2)n\\
  &\phantom{=}+k^2+u^2+uk+vk+wk+uv+uw+vw+k-u+1\\
  B(n)&=-n(n+u)(n+u+v)(n+u+w)\;,
\end{align*}
where $u$, $v$, $w$, and $k$ are nonnegative integer parameters, with
speed of convergence in $C/n^{u+v+w+2k+1}$ for some constant $C$.

An extensive additional search seems to show that there are no additional
sporadic representations. However, we will see in Section \ref{sec:half} that
it is also possible to choose $v$ or $w$ (but not both) in $1/2+\Z$.

\medskip

{\bf Proof:} in the present case, it is possible, not only to prove that
the four-parameter family above does represent $\pi^2$, but also to
explain where it originates.

Consider the series $S=\sum_{n\ge1}1/\prod_{0\le i\le m-1}(n+i)^{e_i}$,
where $e_i=1$ if $0\le i\le m_1-1$ or $m_1+m_2\le i\le m-1$,
and $e_i=2$ if $m_1\le i\le m_1+m_2-1$, where $m_1\ge0$, $m_2\ge1$,
and $m\ge m_1+m_2$. By expanding into partial fractions it is immediate
to show that the sum $S$ satisfies
$S=a\pi^2+b$ with $a$, $b$ rational and $a\ne0$. On the other hand,
applying Euler's transformation of $S$ into a CF, and simplifying by
suitable factors, it is not difficult to show that $S$ (hence $\pi^2$) is
represented by the CF $(A(n),B(n))$, with
\begin{align*}A(n)&=(n-1)(n+m_1-1)+(n+m_1+m_2-1)(n+m-1)\text{\quad and\quad}\\
B(n)&=-n(n+m_1)(n+m_1+m_2-1)(n+m-1)\;.\end{align*}
This gives the three-parameter family above when $k=0$, and the case
of general $k$ is obtained by successive Bauer--Muir accelerations.

\medskip

\subsection{Second Family}

\medskip

We also have the modified series obtained by removing the Euler factor
at $2$: $\pi^2/8=\sum_{n\ge1}1/(2n-1)^2$, giving the
CF $((0,1,8n^2-16n+10),(1,-(2n-1)^4))$ for $\pi^2/8$, that we write in
abbreviated form $$\pi^2=(8n^2-16n+10,-(2n-1)^4)\;.$$
All the work has been done for us in the first family: we simply replace
$n$ by $n-1/2$ and multiply by suitable powers of $2$ to simplify the
formulas, giving the single four-parameter family $\pi^2=(A(n),B(n))$ with
\begin{align*}
  A(n)&=8n^2+4(3u+v+w-4)n\\
  &\phantom{=}+4(k^2+u^2+uk+vk+wk+uv+uw+vw)+4k-10u-2v-2w+10\\
  B(n)&=-(2n-1)(2n-1+2u)(2n-1+2u+2v)(2n-1+2u+2w)\;,\end{align*}
with the same speed of convergence in $C/n^{u+v+w+2k+1}$.

However, inspection shows that while $v$, $w$, and $k$ must still be
nonnegative, we can have $u$ negative as long as $u\ge-k-\min(v,w)$.
Note that here we \emph{cannot} choose $v$ or $w$ in $1/2+\Z$.

\medskip

\section{CFs for $\pi^2$: Convergence in $1/n^P$, $P\in1/2+\Z_{\ge0}$}\label{sec:half}

\medskip

We recall the following formula due to Euler:
$$\sum_{n\ge1}\dfrac{(2z)^{2n}}{n^2\binom{2n}{n}}=2\asin^2(z)\;,$$
which we will use several times. Here, we chooose $z=1$ to obtain
the series
$$\sum_{n\ge1}\dfrac{4^n}{n^2\binom{2n}{n}}=\dfrac{\pi^2}{2}\;,$$
and using Euler's transformation of series into CFs we see that $\pi^2/2$
is given by the CF $\pi^2/2=((0,4n^2-5n+2),(2,-2n^3(2n-1)))$ with convergence
identical to that of the series in $\sqrt{\pi}/n^{1/2}$, which we write
in abbreviated form
$$\pi^2=(4n^2-5n+2,-2n^3(2n-1))\;.$$
If we apply the Bauer--Muir transformation, we check that for any nonnegative
integer $k\ge0$ we also have $\pi^2=(4n^2-5n+2k^2+k+2,-2n^3(2n-1))$.

Here, we could do a new search, but luckily we have already done the work:
it is simply the first four-parameter family found above, but now with the
parameter $w$ being a half integer instead of an integer. We thus replace
$w$ by $w-1/2$ and simplify, so we obtain the four-parameter family
$\pi^2=(A(n),B(n))$ with
\begin{align*}
  A(n)&=4n^2+(6u+2v+2w-5)n+2k^2\\
  &\phantom{=}+(2u+2v+2w+1)k+2u^2+(2v+2w-3)u+v(2w-1)+2\\
  B(n)&=-2n(n+u)(n+u+v)(2n+2u+2w-1)\;,
\end{align*}
where $u$, $v$, and $k$ are nonnegative integer parameters and $w\in\Z$, with
speed of convergence in $C/n^{u+v+w+2k+1/2}$ for some constant $C$.

The \emph{proof} of the validity of this family follows from the one
given in Section \ref{sec:first}.

\medskip

\section{CFs for $\pi^2$: Convergence in $(-1)^n/n^P$}

\medskip

Using Euler's transformation of series into CFs but now applied to
$\pi^2/12=\sum_{n\ge1}(-1)^{n-1}/n^2$, we see that $\pi^2/12$ is given by the
CF $((0,2n-1),(1,n^4))$, with speed of convergence identical to that of the
series, in $(-1)^n/(2n^2)$, which we write in abbreviated form
$$\pi^2=(2n-1,n^4)\;.$$
If we apply the Bauer--Muir transformation, we immediately see that for
any \emph{odd} integer $k\ge1$ we also have $\pi^2=(k(2n-1),n^4)$.
We easily check that this is also true for even integers, so we have found
a one-parameter family. Now we search numerically for other
representations of the form $\pi^2=(A(n),B(n))$, with $A(n)$ of degree $1$ and
$B(n)$ monic of degree $4$, as before assumed to be a product of linear
factors with integral coefficients. Since $(-A(n),B(n))$
is easily seen to be equivalent to $(A(n),B(n))$, we can also
assume that the leading coefficient of $A(n)$ is positive, so we set
$A(n)=a_1n+a_0$ with $a_1\ge1$ and $B(n)=n(n+u_1)(n+u_2)(n+u_3)$ with
$0\le u_1\le u_2\le u_3$. We again find a large number of representations,
and after some cleaning up, we obtain a three-parameter family
and four one-parameter families, and no additional sporadic representations:
\begin{align*}
\pi^2&=((k+1)(2n+2u+2v-1),n(n+u)(n+u+2v)(n+2u+2v))\\
&=((2u^2+8u+7)n+(2u+3)(u^2+3u+1),n(n+u)(n+u+1)(n+2u+4))\\
&=((2u^2+8u+7)n+(2u+3)(u+2)(u+3),n(n+u+3)(n+u+4)(n+2u+4))\\
&=((8u+7)n+(24u^2+24u+3),n(n+2u)(n+4u+1)(n+6u+4))\\
&=((8u+7)n+6(u+1)(4u+3),n(n+2u+3)(n+4u+4)(n+6u+4))\;,\end{align*}
where $u$, $v$, and $k$ are nonnegative integer parameters, with
speeds of convergence $(-1)^nC/n^{2k+2}$, $(-1)^nC/n^{2u^2+8u+7}$,
$(-1)^nC/n^{2u^2+8u+7}$, $(-1)^nC/n^{8u+7}$ and $(-1)^nC/n^{8u+7}$
respectively, for some constants $C$.

\medskip

{\bf Proof:} Although similar to the proof given in Section \ref{sec:first}
the proof of the above is both more complicated and more
interesting since we will need to generalize the construction. Here we
consider the series
$$S=\sum_{n\ge1}(-1)^{n-1}/\prod_{0\le i\le m-1}(n+i)^{e_i}\;,$$
with the same $e_i$ as above. Once again we can simplify the result
given by Euler's transformation of $S$ into a CF, and we obtain that
$S$ is represented by the CF $(A(n),B(n))$ with
\begin{align*}A(n)&=-(n-1)(n+m_1-1)+(n+m_1+m_2-1)(n+m-1)\text{\quad and\quad}\\
B(n)&=n(n+m_1)(n+m_1+m_2-1)(n+m-1)\;.\end{align*}
The problem is that the sum $S$ is not anymore of the form $a\pi^2+b$, but
of the form $a\pi^2+c\log(2)+b$. Thus, if we want a CF representing $\pi^2$
itself, we need $c=0$. It is not difficult to show that this is indeed the
case when $(m_1,m_2,m_3):=(m_1,m_2,m-m_1-m_2)$ is of one of the following
three shapes: $(m_1,m_2,m_3)=(m_1,2u+1,m_1)$ or $(2u,2u+2,2u+3)$, or
$(2u+3,2u+2,2u)$. It should not be difficult to prove that there are no
other cases, but I have not done so. In any event, the first shape
followed by Bauer--Muir accelerations leads to the three-parameter family,
and the other two shapes lead to the last two one-parameter families.

The first two one-parameter families do \emph{not} come from this
construction, but from a more general one. We first note the following
result, whose easy proof is left to the reader

\begin{proposition}\label{prop1}
  Let $c(x)\in\Q[x]$ be a polynomial with $\deg(c)\ge2$,
  let $z=\pm1$, let $P\in\Q[x]$ be an irreducible polynomial which divides
  $c(x)P(x+1)+zc(x-1)P(x-1)$, let $R(x)=(c(x)P(x+1)+zc(x-1)P(x-1))/P(x)$
  and $g(x)=\gcd(c(x),c(x-1))$, and set
  $$S=\sum_{n\ge1}\dfrac{z^n}{c(n)P(n)P(n+1)}\;.$$
  Assume that $c(x)$, $P(x)$, and $P(x+1)$ are pairwise coprime.
  We have the continued fraction
  $S=[[0,R(n)/g(n)],[r_1,-zc(n)^2/(g(n)g(n+1))]]$ and the identity
  $S=r_2+\sum_{n\ge1}z^nd(n)/c(n)$ for some rational numbers $r_1$ and
  $r_2$ and polynomial $d$ with $\deg(d)<\deg(c)$.
\end{proposition}  

We also need the following conjecture:

\begin{conjecture}
  Let $u$ be a nonnegative integers. Up to a multiplicative constant there
  exists a unique polynomial $P$ such that
  \begin{align*}(x+2u+4)(x+u+1)&P(x+1)-(x-1)(x+u-1)P(x-1)\\
    &=((2u^2+8u+7)x+(2u+3)(u^2+3u+1))P(x)\;,\end{align*}
  and we have $\deg(P)=u(u+3)$. The same is true if we require
  \begin{align*}(x+2u+4)(x+u+4)&P(x+1)-(x-1)(x+u+2)P(x-1)\\
    &=((2u^2+8u+7)x+(2u+3)(u+2)(u+3))P(x)\;.\end{align*}
\end{conjecture}

Assuming this conjecture, in Proposition \ref{prop1} we choose $z=-1$,
$c(x)=\prod_{0\le i\le 2u+4}(n+i)^{e_i}$ as before, with $m_1=u$, $m_2=2$, and
$m_3=u+3$, so $m=2u+5$, and $P$ the polynomial satisfying the first identity
above. We check that $g(x)=(x+u)\prod_{0\le i\le 2u+3}(x+i)$, hence
$$-zc(n)^2/(g(n)g(n+1))=n(n+u)(n+u+1)(n+2u+4)\;,$$
and because of the identity for $P$,
$$R(n)/g(n)=((2u^2+8u+7)n+(2u+3)(u^2+3u+1)\;.$$
Finally, the identity for
$S$ given by Proposition \ref{prop1} implies that $S$ is a linear combination
of $\pi^2$, $\log(2)$, and $1$. Similarly for the second identity.
Thus, the other two one-parameter families follow from the above conjecture
together with the additional conjecture that $\log(2)$ does not occur
in the linear combination for $S$.

\smallskip

{\bf Remark:} Concerning the above conjecture, Peter Mueller from the
MathOverflow forum proves the following:
\begin{enumerate}
\item The two statements are equivalent.
\item If $P$ exists it is indeed unique up to a multiplicative constant,
  and has degree $u(u+3)$.
\end{enumerate}

\smallskip

It is possible that there exist other, possibly more complicated,
one-parameter families.

\section{CFs for $\pi^2$: Convergence in $1/((4/3)^nn^P)$}\label{sec:43}

Here we use the $\asin^2(z)$ series with $z=\sqrt{3}/2$, giving the series
$$\sum_{n\ge1}\dfrac{3^n}{n^2\binom{2n}{n}}=\dfrac{2\pi^2}{9}\;,$$
and applying Euler's transformation shows that $2\pi^2/9$ is given by the
CF $((0,7n^2-8n+3),(3,-6n^3(2n-1)))$ with speed of convergence identical
to that of the series, in $3\sqrt{\pi}/((4/3)^nn^{3/2})$, and as before
we can write in abbreviated form
$$\pi^2=(7n^2-8n+3,-6n^3(2n-1))\;.$$
Applying Bauer--Muir does not give anything simple. So as before we search
for representations $\pi^2=(A(n),B(n))$ with $A(n)=7n^2+a_1n+a_0$ and
$B(n)=-6n(n+u_1)(n+u_2)(2n-1+2u_3)$. Note that we could try other types
for $B(n)$ with the same leading term $-12n^4$ and completely factored,
for instance $B(n)=-n(2n+u_1)(2n+u_2)(3n+u_3)$, but I have not found any
new representations in this way.

Here, the situation is less favorable since I have not been able to
classify completely the representations. I have found four two-parameter
families:

\begin{align*}
  \pi^2&=(7n^2+(2u+3v-6)n-(5u^2+15vu+18v^2-2),\\
  &\phantom{=}-6n(n+u)(n+u+3v)(2n+1-4u-6v))\\
  &=(7n^2+(2u+3v-8)n-(5u^2+(15v+2)u+3(3v-1)(2v+1)),\\
  &\phantom{=}-6n(n+u)(n+u+3v)(2n-1-4u-6v))\\
  &=(7n^2+(15u+3v-8)n-3(2v+1)(3u+3v-1),\\
  &\phantom{=}-6n(n+u)(n+3u+3v)(2n-1-6v))\\
  &=(7n^2+(15u+3v-3)n-2(3v+4)(3u+3v+2),\\
  &\phantom{=}-6n(n+u)(n+3u+3v+3)(2n-5-6v))\\
\end{align*}
with all parameters nonnegative, plus a number of apparently sporadic cases:

\begin{align*}
  \pi^2
  &=(7n^2-16n+3,-6n^3(2n-5))\\
  &=(7n^2-19n-6,-6n^2(n+1)(2n-11))\\
  &=(7n^2+17n,-6n(n+3)^2(2n+1))\\
  &=(7n^2+14n-36,-6n(n+3)(n+5)(2n-7))\\
  &=(7n^2+54n+5,-6n(n+6)(n+10)(2n+1))\\
  &=(7n^2+31n-60,-6n(n+1)(n+10)(2n-3)) \;,\end{align*}


all with convergence speeds $C/((4/3)^nn^P)$ for varying half-integral $P$.

\smallskip

I have not tried to prove the validity of the above CFs.

The upshot of this case is that first, that I have probably not found all
possible cases, and second, that I have not seen the patterns leading
either to additional two-parameter families, or even possibly to a
three-parameter family.

Note, however, that in Section \ref{sec:44} we will meet the same
\emph{exact} situation (four two-parameter families plus six sporadic cases),
which makes me believe that there is a relation between the CFs of
that section and those of the present section.

\section{CFs for $\pi^2$: Convergence in $1/(2^nn^P)$}

Here we use the $\asin^2(z)$ series with $z=\sqrt{2}/2$, giving the series
$$\sum_{n\ge1}\dfrac{2^n}{n^2\binom{2n}{n}}=\dfrac{\pi^2}{8}\;,$$
and applying Euler's transformation shows that $\pi^2/8$ is given by the
CF $((0,3n^2-3n+1),(1,-n^3(2n-1)))$ with speed of convergence identical
to that of the series, in $\sqrt{\pi}/(2^nn^{3/2})$, and as before
we can write in abbreviated form
$$\pi^2=(3n^2-3n+1,-n^3(2n-1))\;.$$

Here, applying Bauer--Muir gives the one-parameter family
$\pi^2=(3n^2-3n+1,-n^3(2n-1-2k))$, and we search
for representations $\pi^2=(A(n),B(n))$ with $A(n)=3n^2+a_1n+a_0$ and
$B(n)=-n(n+u_1)(n+u_2)(2n-1+2u_3)$. After some work, we find
the three-parameter family:
\begin{align*}\pi^2&=(3n^2+(8u+4v-3)n+(2u-1)(2u+2v-1),\\
  &\phantom{=}-n(n+2u)(n+2u+2v)(2n+2u-2k-1))\;,\end{align*}
where all the parameters are nonnegative, with speed of convergence
$$C/(2^nn^{u+2v+3/2+3k})$$ for some constant $C$.

An extensive additional search seems to show that there are no additional
sporadic representations.

It is easy to give a \emph{proof} of the validity of the above family,
which we leave to the reader.

\section{CFs for $\pi^2$: Convergence in $1/(4^nn^P)$}\label{sec:44}

Here we use the $\asin^2(z)$ series with $z=1/2$, giving the series
$$\sum_{n\ge1}\dfrac{1}{n^2\binom{2n}{n}}=\dfrac{\pi^2}{18}\;,$$
and applying Euler's transformation shows that $\pi^2/18$ is given by the
CF $((0,5n^2-4n+1),(1,-2n^3(2n-1)))$ with speed of convergence identical
to that of the series, in $\sqrt{\pi}/3/(4^nn^{3/2})$, and as before
we can write in abbreviated form
$$\pi^2=(5n^2-4n+1,-2n^3(2n-1))\;.$$

Applying Bauer--Muir leads to an explicit complicated one-parameter family
which will be a special case of the ones below, and the situation is
similar to that of Section \ref{sec:43}: I have found four two-parameter
families:

\begin{align*}
  \pi^2&=(5n^2+(14u+21v-4)n+(3u+3v-1)(3u+6v-1),\\
  &\phantom{=}-2n(n+u)(n+u+3v)(2n-1))\\
  &=(5n^2+(14u+21v-6)n+9u^2+(27v-8)u+2(3v-1)^2,\\
  &\phantom{=}-2n(n+u)(n+u+3v)(2n+1))\\
  &=(5n^2+(17u+21v+15)n+2(2u+3v+2)(3u+3v+2),\\
  &\phantom{=}-2n(n+u)(n+3u+3v+3)(2n+2u+1))\\
  &=(5n^2+(17u+21v-4)n+(3u+3v-1)(4u+6v-1),\\
  &\phantom{=}-2n(n+u)(n+3u+3v)(2n+2u-1))\;,
  \end{align*}
with all parameters nonnegative, plus a number of apparently sporadic cases:

\begin{align*}
  \pi^2&=(5n^2+4n+1,-2n^3(2n+3))\\
  &=(5n^2+15n+6,-2n^2(n+1)(2n-1))\\
  &=(5n^2+19n+16,-2n(n+3)^2(2n+3))\\
  &=(5n^2+78n+261,-2n(n+3)(n+7)(2n-1))\\
  &=(5n^2+38n+68,-2n(n+3)(n+5)(2n+1))\\
  &=(5n^2+45n+60,-2n(n+1)(n+10)(2n+3))\;,\end{align*}
all with convergence speeds $C/(4^nn^P)$ for varying half-integral $P$.

As in Section \ref{sec:43}, I have not tried to prove the validity of these
CFs, I have also probably not found all possible cases, and I have not
seen the patterns leading either to additional one- or two-parameter families,
or even possibly to a three-parameter family.

\smallskip

Nonetheless, the reader has certainly noticed that there is a striking
resemblance between the families of this section and those of Section
\ref{sec:43}, and this should probably be investigated.

\section{CFs for $\pi^2$: Convergence in $(-1)^n/(8^nn^P)$}

Let $$F_n=\sum_{0\le k\le n}\binom{n}{k}^3$$
be the $n$th \emph{Franel number}. From the recursion found by Franel
for these numbers and the computation of so-called Ap\'ery limits
(see for instance \cite{Gor}),
we have the CF $\pi^2/3=((0,7n^2-7n+2),(8,8n^4))$ with speed of
convergence $C/(-8)^n$ and whose $n$th denominator
$q(n)$ is given by $q(n)=n!^2F_n$, giving the formula
$$\dfrac{\pi^2}{3}=\sum_{n\ge1}(-1)^{n-1}\dfrac{8^n}{n^2F_nF_{n-1}}\;.$$
As usual we write the CF in abbreviated form
$$\pi^2=(7n^2-7n+2,8n^4)\;.$$
Bauer--Muir acceleration gives something complicated, so as usual we
search. The situation is even worse here: I have found experimentally a
large number of representations, for which I was able to extract three
two-parameter families and a small number of one-parameter families,
as follows:

\begin{align*}
\pi^2&=(7n^2+((48-20v)u-7)n+2(4u-1)((8-5v)u-1),\\
     &\phantom{=((}8n(n+uv)(n+u(v+2))(n+4u))\\
     &=(7n^2+(24v-20u-7)n+(4v-2)(4v-5u-1),\\
     &\phantom{=((}8n(n+u)(n+2v)(n+u+v))\\
     &=(7n^2+(12u+20v-7)n+(4v-2)(3u+3v-1),\\
     &\phantom{=((}8n(n+2v)(n+3u+2v)(n+3u+3v))\\
     &=(7n^2+(30u-7)n+(3u-1)(9u-2),8n^2(n+3u)^2)\\
     &=(7n^2+(34u-7)n+27u^2-17u+2,8n^2(n+u)^2)\\
     &=(7n^2+(22u-7)n+(3u-1)(5u-2),8n(n+u)^2(n+3u))\\
     &=(7n^2+(58u-7)n+87u^2-29u+2,8n(n+u)^2(n+3u))\\
     &=(7n^2+(48u-23)n+2(4u+1)(8u-7),8n(n+2)(n+2u+3)(n+4u+2))\\
     &=(7n^2+(34u+17)n+3(3u+1)(3u+2),8n(n+1)(n+u)(n+u+2))\;,
\end{align*}
where $u$ and $v$ are nonnegative integers, except in the third two-parameter
representation where one can have $u\ge-\lfloor 2v/3\rfloor$, all with
convergence speeds $C/((-8)^nn^P)$ for varying integral $P$.

I have found many (for now more than $60$) additional individual CFs, from
which one can probably extract more one-parameter families. Here is a small
sample:

\begin{align*}
  \pi^2&=(7n^2+26n+8,8n^3(n+3))\\
  &=(7n^2+80n+105,8n^2(n+1)(n+2))\\
  &=(7n^2+24n+9,8n^2(n+1)(n+4))\\
  &=(7n^2+65n+70,8n^2(n+1)(n+8))\\
  &=(7n^2-15n-12,8n^2(n+4)^2)\\
  &=(7n^2+57n+96,8n^2(n+4)^2)\;.\end{align*}

\section{CFs for $\pi^2$: Convergence in $1/(-((1+\sqrt{5})/2)^{10})^n$}

This is the famous CF of Ap\'ery:

$$\pi^2=((0,11n^2-11n+3),(30,n^4))$$

with convergence in $24\pi^2/(1+\sqrt{2})^{10n+5}$.

This CF seems to be unique of its kind, i.e. of the form
$\pi^2=(11n^2+a_1n+a_0,n(n+u_1)(n+u_2)(n+u_3))$.

\section{CFs for $\pi^2$: Other Speeds of Convergence}

I have not found any convergent CF of the form
$\pi^2=(A(n),B(n))$ with $\deg(A)\le 2$, $\deg(B)\le 4$ and $A$ and $B$
with integral coefficients other than the
families mentioned above, up to trivial equivalence such as changing
$n$ into $n+j$ or multiplying $(A,B)$ by $(c,c^2)$ for a nonzero
constant $c$. Of course, there exist infinitely many other
CFs for $\pi^2$, but of larger bidegrees.

\section{CFs for $\z(3)$: Convergence in $1/n^P$}

We now study CFs for $\z(3)$. Similarly to $\pi^2$, as far as the author
is aware, up to equivalence there exist exactly eight types of convergence of
polynomial continued fractions for $\z(3)$ with integral coefficients and
bidegree $\le(3,6)$:
three families with convergence in $1/n^P$ with $P\in\Z$, a single CF
with convergence in $(-1)^n/n^P$, a single CF with convergence in
$1/(4^nn^P)$, a family of CFs with convergence in $1/((-4)^nn^P)$,
a single CF with convergence in $1/E^n$ with $(1+\sqrt{2})^4$, and finally
the unique CF of Ap\'ery with convergence in $1/E^n$ for $E=(1+\sqrt{2})^8$.
We consider each of these cases in turn.

\medskip

\subsection{First Family}\label{sec:firstz3}

\medskip

The natural series giving $\z(3)$ is $\z(3)=\sum_{n\ge1}1/n^3$, so Euler's
transformation of series gives the CF $\z(3)=((0,(2n-1)(n^2-n+1)),(1,-n^6))$,
abbreviated as usual as $$\z(3)=((2n-1)(n^2-n+1),-n^6)\;.$$ Bauer--Muir
acceleration leads to the family $\z(3)=((2n-1)(n^2-n+2k^2+2k+1),-n^6)$
for nonnegative integer $k$.

Instead of doing a search, we generalize the proof given in Section
\ref{sec:first}: using the same notation, we choose $e_i=1$, $2$, $3$, $2$, or
$1$ according to $0\le i<m_1$, $m_1\le i<m_1+m_2$, $m_1+m_2\le i<m_1+m_2+m_3$,
$m_1+m_2+m_3\le i<m_1+m_2+m_3+m_4$, or $m_1+m_2+m_3+m_4\le i<m$ respectively,
and this guarantees that the CF $(A(n),B(n))$ will have bidegree $\le(3,6)$,
with
\begin{align*}
    A(n)&=(n+m_1+m_2+m_3-1)(n+m_1+m_2+m_3+m_4-1)(n+m-1)\\
  &\phantom{=}+(n-1)(n+m_1-1)(n+m_1+m_2-1)\text{\quad and}\\
    B(n)&=-n(n+m_1)(n+m_1+m_2)(n+m_1+m_2+m_3-1)\\
    &\phantom{=}(n+m_1+m_2+m_3+m_4-1)(n+m-1)\;.
\end{align*}

However the sum $S$ will be a linear combination of $1$, $\pi^2$, and $\z(3)$,
and we do not want any $\pi^2$. For this, a \emph{sufficient} condition
is that $m_3$ is odd and $m_4=m_2$ and $m_5=m_1$, giving the parametric
family $(A,B)$ with
\begin{align*}
  A(n)&=(n+u+v+2w)(n+u+2v+2w)(n+2u+2v+2w)\\
  &\phantom{=}+(n-1)(n+u-1)(n+u+v-1)\text{\quad and}\\
  B(n)&=-n(n+u)(n+u+v)(n+u+v+2w)(n+u+2v+2w)(n+2u+2v+2w)\;,\end{align*}
and note that
\begin{align*}A(n)&=(2n+2u+2v+2w-1)(n^2+(2u+2v+2w-1)n\\
  &\phantom{=}+(u^2+(3v+4w)u+2v^2+(6w+1)v+4w^2+2w+1))\;.\end{align*}
Applying Bauer--Muir, we thus find the four-parameter family
$\z(3)=(A(n),B(n))$ with
\begin{align*}
  A(n)&=(2n+2(u+v+w)-1)(n^2+(2(u+v+w)-1)n+u^2+2v^2+4w^2\\
   &\phantom{=}+3uv+4uw+6vw+v+2w+1+2k(u+2v+3w)+2k^2+2k),\\
 B(n)&=-n(n+u)(n+u+v)(n+u+v+2w)(n+u+2v+2w)(n+2u+2v+2w)\;,\end{align*}
where $u$, $v$, $w$, and $k$ are nonnegative integer parameters, with
speed of convergence in $C/n^{2u+4v+6w+4k+2}$ for some constant $C$.

However, the sufficient condition (symmetry of $(m_1,m_2,m_3,m_4,m_5)$)
above is not necessary: it is easy to find as many nonsymmetrical examples
as we like, the simplest being $(m_1,m_2,m_3,m_4,m_5)=(0,0,2,2,0)$.
I have not yet been able to give a reasonable classification of these
additional examples, here is a small sample:
\begin{align*}
\z(3)&=(2n^3+4n^2+18n+8,-n^3(n+1)(n+3)^2)\\
     &=(2n^3+8n^2+26n+28,-n^2(n+2)(n+3)^3)\\
     &=(2n^3+11n^2+53n+60,-n^2(n+1)(n+2)(n+5)(n+6))\\
     &=(2n^3+17n^2+78n+119,-n(n+2)^2(n+3)(n+5)(n+8))\\
     &=(2n^3+19n^2+93n+180,-n(n+1)(n+4)(n+5)(n+6)^2)\\
     &=(2n^3+20n^2+163n+299,-n^3(n+3)(n+10)^2)\\
     &=(2n^3+22n^2+156n+360,-n(n+1)(n+2)(n+5)(n+8)(n+9))\\
     &=(2n^3+23n^2+201n+500,-n^2(n+1)(n+5)(n+10)^2)\;.
\end{align*}
\medskip

\subsection{Second Family}

We also have the simple modification obtained by removing the Euler
factor at $2$: $\z(3)=(8/7)\sum_{n\ge1}1/(2n-1)^3$, so Euler's
transformation gives the CF $\z(3)=((0,1,4(n-1)(4n^2-8n+7)),(8/7,-(2n-1)^6))$,
abbreviated as usual as $$\z(3)=(4(n-1)(4n^2-8n+7),-(2n-1)^6)\;.$$ Bauer--Muir
acceleration leads to the family $\z(3)=(4(n-1)(4n^2-8n+8k^2+8k+7),-(2n-1)^6)$
for nonnegative integer $k$.
As in the case of $\pi^2$, to find parametric families we simply replace
$n$ by $n-1/2$ and multiply by suitable powers of $2$ in the first family
found above, giving the four-parameter family $\z(3)=(A(n),B(n))$ with

\begin{align*}
  A(n)&=4(n+u+v+w-1)(4n^2+8(u+v+w-1)n+4u^2+8v^2+16w^2\\
   &\phantom{=}+12uv+16uw+24vw-4u+4w+8k(u+2v+3w)+8k^2+8k+7),\\
  B(n)&=-(2n-1)(2n+u-1)(2n+u+v-1)(2n+u+v+2w-1)\cdot\\
  &\phantom{=-}\cdot(2n+u+2v+2w-1)(2n+2u+2v+2w-1)\;,\end{align*}
where $u$, $v$, $w$, and $k$ are nonnegative integer parameters, with
speed of convergence in $C/n^{2u+4v+6w+4k+2}$ for some constant $C$,
and of course as in the first family, one also has the numerous additional
apparently sporadic cases, such as
$$\z(3)=(16n^3+8n^2+124n-2,-(2n-1)^3(2n+1)(2n+5)^2)\;.$$

\subsection{Third Family}

From work of Y.~Yang \cite{Yan}, one finds the CF
$$\z(3)=((0,(2n-1)(2n^2-2n+1)),(2/7,-4n^6))$$ with extremely slow convergence
in $C/\log(n)$, which as usual we abbreviate to
$\z(3)=((2n-1)(2n^2-2n+1),-4n^6)$. Modifying his construction and applying
Bauer--Muir, we find the two-parameter family $\z(3)=(A(n),B(n))$ with
\begin{align*}
  A(n)&=(2n+2u-1)(2n^2+(4u-2)n+(2k+2)u+k^2+4k+5),\\
  B(n)&=-4n^3(n+2u)^3\;,
\end{align*}
with $u$ and $k$ nonnegative, with speed of convergence in $C/n^{2u+2k+4}$.

\smallskip

It is interesting to compare this with the very similar two-parameter family
$\z(3)=(A(n),B(n))$ obtained by setting $u=v=0$ and replacing $w$ by $u$
in Section \ref{sec:firstz3}:
\begin{align*}
  A(n)&=(2n+2u-1)(2n^2+(4u-2)n+8u^2+(12k+4)u+4k^2+4k+2),\\
  B(n)&=-4n^3(n+2u)^3\;,
\end{align*}
with $u$ and $k$ nonnegative, with speed of convergence in $C/n^{6u+4k+2}$.

\bigskip

\section{CFs for $\z(3)$: Convergence in $(-1)^n/n^P$}

Using Euler's transformation of series into CFs but now applied to
$(3/4)\z(3)=\sum_{n\ge1}(-1)^{n-1}/n^3$, we see that $\z(3)$ is given by the
CF $((0,3n^2-3n+1),(4/3,n^6))$, with speed of convergence identical to that of
the series, in $(-1)^n2/(3n^3)$, which we write in abbreviated form
$$\z(3)=(3n^2-3n+1,n^6)\;.$$
There is apparently no other CF of this type with convergence in $(-1)^n/n^P$.
The main reason is probably that, if we try to apply similar constructions
to those that we used for $\sum_{n\ge1}(-1)^{n-1}/n^2$, we are now going
to obtain linear combinations of $1$, $\z(3)$, and the two additional
``parasitic'' terms $\z(2)$ and $\log(2)$, instead of simply $\log(2)$ alone.

\section{CFs for $\z(3)$: Convergence in $1/(4^nn^P)$}

Let $$D_n=\sum_{0\le k\le n}\binom{n}{k}^2\binom{2k}{k}\binom{2n-2k}{n-k}$$
be the $n$th so-called \emph{Domb number}. From the recursion that one can find
for these numbers and the computation of Ap\'ery limits
(see for instance \cite{Gor}), we have the CF
$\z(3)=((0,(2n-1)(5n^2-5n+2)),(12/7,-16n^6))$, whose $n$th denominator
$q(n)$ is given by $q(n)=n!^3D_n/2^n$, giving the formula
$$\dfrac{56\z(3)}{3}=\sum_{n\ge1}\dfrac{64^n}{n^3D_nD_{n-1}}\;.$$
As usual we write the CF in abbreviated form
$$\z(3)=((2n-1)(5n^2-5n+2),-16n^6)\;.$$

Bauer--Muir acceleration gives something complicated, and as above I have not
found any other CF of this type with convergence in $1/(4^nn^P)$.

\section{CFs for $\z(3)$: Convergence in $(-1)^n/(4^nn^P)$}

From the famous series mentioned by Ap\'ery but known before him
$$\z(3)=\dfrac{5}{2}\sum_{n\ge1}\dfrac{(-1)^{n-1}}{n^3\binom{2n}{n}}\;,$$
one deduces the CF $\z(3)=((0,3n^3+n^2-3n+1),(5/2,2n^5(2n-1)))$ with the
same speed of convergence as the series in $\sqrt{\pi}/2/((-4)^nn^{5/2})$,
which we write as
$$\z(3)=(3n^3+n^2-3n+1,2n^5(2n-1))\;.$$
Here, a variant of Bauer--Muir gives the family
$$\z(3)=(3n^3+(11k+1)n^2-(11k+3)n+3k+1,n^4(2n-k)(2n-k-1))$$
with convergence in $C/((-4)^nn^{5k+5/2})$. A search apparently does not
give any other CF.

\section{CFs for $\z(3)$: Convergence in $1/(1+\sqrt{2})^{4n}$}

Let $$C_n=\sum_{n/2\le k\le n}\binom{n}{k}^2\binom{2k}{n}^2$$
(these numbers do not seem to have a name).
From the recursion that one can find
for these numbers and the computation of Ap\'ery limits
(see for instance \cite{Gor}), we have the CF
$\z(3)=((0,(2n-1)(3n^2-3n+1)),(8/7,-n^6))$ with speed of convergence
$4\pi^3/7/(1+\sqrt{2})^{4n+2}$, whose $n$th denominator
$q(n)$ is given by $q(n)=n!^3C_n/4^n$, giving the formula
$$\dfrac{7\z(3)}{2}=\sum_{n\ge1}\dfrac{16^n}{n^3C_nC_{n-1}}\;.$$
As usual we write the CF in abbreviated form
$$\z(3)=((2n-1)(3n^2-3n+1),-n^6)\;.$$

Since the exponent of the speed of convergence is irrational, Bauer--Muir
acceleration would give irrational coefficients so is not applicable,
and I have not been able to find other CFs of this type.

\section{CFs for $\z(3)$: Convergence in $1/(1+\sqrt{2})^{8n}$}

This is the famous CF of Ap\'ery:

$$\z(3)=((0,(2n-1)(17n^2-17n+5)),(6,-n^6))$$
with convergence in $4\pi^3/(1+\sqrt{2})^{8n+4}$.

This CF seems to be unique of its kind, i.e. of the form
$\z(3)=(34n^3+a_2n^2+a_1n+a_0,-n(n+u_1)(n+u_2)(n+u_3)(n+u_4)(n+u_5)(n+u_6))$.

\section{CFs for $\z(3)$: Other Speeds of Convergence}

I have not found any convergent CF of the form
$\z(3)=(A(n),B(n))$ with $\deg(A)\le 3$, $\deg(B)\le 6$ and $A$ and $B$
with integral coefficients other than the
families mentioned above, up to trivial equivalence such as changing
$n$ into $n+j$ or multiplying $(A,B)$ by $(c,c^2)$ for a nonzero
constant $c$. Of course, there exist infinitely many other
CFs for $\z(3)$, but of larger bidegrees.

\section{Other Periods}

All the preceding sections dealt with CFs for $\pi^2$ and $\z(3)$.
A similar study can be done for other periods. I briefly summarize what I
have found, restricting as usual to CFs of bidegrees $\le(k,2k)$ for periods
of degree $k$.

\smallskip

\subsection{Periods of Degree $1$}

\smallskip

As just mentioned, we restrict to bidegrees $\le(1,2)$.

$\bullet$ $z=\pi$: one finds parametric families of CFs with convergence in
$(-1)^n/n^P$, $1/(2^nn^P)$, and $(-1)^n/(1+\sqrt{2})^{2n}$.

$\bullet$ $z=\pi/\sqrt{3}$: one finds parametric families of CFs with
convergence in $1/((4/3)^nn^P)$, $(-1)^n/(3^nn^P)$, $1/(4^nn^P)$, and
$(-1)^n/(2+\sqrt{3})^{2n}$.

$\bullet$ $z=\log(2)$: one finds parametric families of CFs with convergence
in $(-1)^n/n^P$, $1/(2^nn^P)$, $(-1)^n/(8^nn^P)$, $1/(9^nn^P)$, and
$1/(1+\sqrt{2})^{4n}$.

Note that one can also study CFs for $z=\log(m)$ for other rational
values of $m$.

$\bullet$ $z=\sqrt{2}\log(1+\sqrt{2})$: one finds parametric families of
CFs with convergence in $1/(1+\sqrt{2})^{2n}$.

\smallskip

\subsection{Periods of Degree $2$}

\smallskip

We restrict to bidegrees $\le(2,4)$.

$\bullet$ $z=\pi^2$: see the extensive lists given in the preceding sections.

$\bullet$ $z=L(\chi_{-3},2)$: one finds parametric families of CFs with
convergence $1/((9/8)^nn^P)$ and $1/(9^nn^P)$.

$\bullet$ $z=L(\chi_{-4},2)=G$, Catalan's constant: one finds parametric
families of CFs with convergence $1/n^P$, $(-1)^n/n^P$, and $1/(2^nn^P)$.

\smallskip

\subsection{Periods of Degree $3$}

We restrict to bidegrees $\le(3,6)$.

$\bullet$ $z=\z(3)$: see the extensive lists given in the preceding
sections.

$\bullet$ For other periods of degree $3$ such as $\pi^3$ and $\pi^3/\sqrt{3}$,
either the CFs are isolated, or are one-parameter families coming from
Bauer--Muir acceleration.

\medskip

I would evidently welcome any additional data and comments
on the subject of this paper.

\end{document}